\documentclass[twoside,10pt]{article}

\usepackage{psfrag}
\usepackage{graphicx}
\usepackage{amssymb}
\usepackage{latexsym}

\newtheorem{theorem}{Theorem}
\newtheorem{proposition}{Proposition}
\newtheorem{definition}{Definition}
\newtheorem{lemma}{Lemma} 
\newtheorem{corollary}{Corollary}

\begin{document}

\title{Non Standard Metric Products}
\maketitle

\vspace{0.2cm}

\begin{center}
{\large Andreas Bernig \footnote{supported by Deutsche Forschungsgemeinschaft}, \hspace{0.3cm} {\large Thomas Foertsch 
\footnote{supported by SNF Grant 21 - 589 38.99},  $\hspace{0.3cm}$ Viktor Schroeder}
}
\end{center}

\vspace{0.5cm}


\begin{abstract}
We consider a fairly general class of natural non standard metric products and classify those amongst them,
which yield a product of certain type (for instance an inner metric space) for all possible 
choices of factors of this type (inner metric spaces). We further prove the additivity of the Minkowski
rank for a large class of metric products.
\end{abstract}


\section{Introduction}

Given a finite number $(X_i,d_i)$, $i=1,...,n$, of metric
spaces there are different possibilities to define a metric $d$ on the product
${\Pi}_{i=1}^{n}X_i$. The standard choice is of course the Euclidean product metric
\begin{displaymath}
d\Big( (x_1,...,x_n),(y_1,...,y_n)\Big) \; = \; 
\Big( \sum\limits_{i=1}^{n} d_i^2(x_i,y_i){\Big)}^{\frac{1}{2}}.
\end{displaymath}

A generalization of this construction is given by warped products (\cite{ab},\cite{chen}), which 
have been proven to be useful for the construction of Hadamard spaces. 
In this paper we consider another generalization, namely the following class of product metrics: \\
Let $(X_i,d_i)$, $i=1,...,n$, be metric spaces and denote the product set by
$X={\Pi}_{i=1}^n X_i$. It is natural to define a metric product  $d$ on $X$ of the form $d=d_{\Phi}$,
\begin{displaymath}
d_{\Phi} \Big( (x_1,...,x_n),(y_1,...,y_n) \Big) \; = \; 
\Phi \Big( d_1(x_1,y_1),...,d_n(x_n,y_n) \Big) ,
\end{displaymath}
where $\Phi :Q^n\longrightarrow [0,\infty )$ is a function defined on the quadrant 
$Q^n=[0,\infty )^n$. \\
Note that for Banach spaces a class of product spaces of this type appears in \cite{cl}. \\ 

The function $\Phi$ has to satisfy certain natural conditions ((A) and (B) in Lemma \ref{metric}) 
in order that $d_{\Phi}$ is a metric.
These conditions still allow strange metrics on the product (even
the trivial product when $n=1$). In particular, $\Phi$ does not have to be continuous. \\

However, once we require for example that $\Phi : Q^n\longrightarrow [0,\infty )$ yields an inner metric space 
$(X,d_{\Phi})$
for all possible choices of inner metric spaces $(X_i,d_i)$, the conditions on $\Phi$ become very rigid.
In fact, those conditions imply that $\Phi$ now has to be continuous. \\
In order to state the corresponding theorem we consider the function
\begin{displaymath}
\Psi : \mathbb{R}^n \longrightarrow [0,\infty ), \hspace{1,5cm}
\Psi \Big( \sum\limits_{i=1}^n \; x_i \, e_i\Big) \; := \; 
\Phi \Big( \sum\limits_{i=1}^n \; |x_i| \, e_i \Big)
\end{displaymath}
and say that $\Phi$ is induced by a norm if and only if the function $\Psi$  is a norm. \\

We say that $\Phi$ preserves length spaces if the product of length spaces, endowed with the metric $d_\Phi$, again is a length space (and correspondingly for other types of metric spaces). 

\begin{theorem} \label{theo-character}
Let $\Phi :Q^n\longrightarrow [0,\infty )$ be a function satisfying conditions (A) and (B) of Lemma \ref{metric}. Then $\Phi$ preserves
\begin{description}
\item[i)] length spaces, 
\item[ii)] geodesic spaces,
\item[iii)] uniquely geodesic spaces,
\item[iv)] metric spaces of non-positive Busemann curvature,
\item[v)] metric spaces of curvature bounded from above (or below),
\end{description}
if and only if
\begin{description}
\item[i),ii)] $\Phi$ is induced by a norm.
\item[iii),iv)] $\Phi$ is induced by a norm with strictly convex norm ball.
\item[v)] $\Phi$ is induced by a scalar product.
\end{description}
\end{theorem}

In Section \ref{sec-minkowski} we consider the behaviour of the Minkowski rank under non standard metric products. 
The Minkowski rank, $rank_M(X,d)$, of a metric space $(X,d)$ is the supremum of the dimensions of normed vector spaces 
isometrically embedded into $X$. We generalize results of \cite{foesch}
to those metric products $(X,d_{\Phi})$ where $\Phi$ is induced by a norm with strictly convex
unit ball:
\begin{theorem} \label{generalized-minkadd}
Let $\Phi : Q^n\longrightarrow [0,\infty )$ be a function such that
$\Psi : \mathbb{R}^n \longrightarrow [0,\infty )$, defined as above,
is a norm with a strictly convex norm ball. Let $(X_i,d_i)$, $i=1,...,n$, be metric spaces
and $X={\Pi}_{i=1}^nX_i$. Then
\begin{displaymath}
rank_{M}\Big( X,d_{\Phi}\Big) \; = \; \sum\limits_{i=1}^n \; rank_{M} \Big( X_i,d_i\Big) .
\end{displaymath} 
\end{theorem}

This shows one advantage of the Minkowski rank over the Euclidean rank, since the latter one is not additive, even with respect to the standard product (\cite{foesch}). \\

A geodesic metric space $(X,d)$ is called convex, if for every pair $c_1:[a_1,b_1]\longrightarrow X$ and
$c_2:[a_2,b_2]\longrightarrow X$ of constant speed geodesics the function 
$d\circ (c_1,c_2):[a_1,b_1]\times :[a_2,b_2] \longrightarrow \mathbb{R}$ is convex. \\
Kleiner (\cite{kl},Theorem D) proved that for locally compact, convex, geodesic metric space $(X,d)$ with cocompactly acting
isometry group the Minkowski rank coincides with a number of other invariants, one of which is the quasi-Euclidean rank. 
The quasi-Euclidean rank, $rank_{qE}(X,d)$, is defined as  
\begin{displaymath}
\sup\{ k \; | \; \mbox{there is a quasi-isometric embedding} \; f:\mathbb{E}^k\longrightarrow X \} .
\end{displaymath}
As a corollary of our results, we
show that these invariants are additive under certain products. In particular we obtain the
\begin{corollary} \label{quasiadd}
Let $\Phi$ be as in Theorem \ref{generalized-minkadd} and $(X_i,d_i)$, $i=1,...,n$, be locally compact convex metric spaces with cocompactly acting isometry group. 
Then for the quasi-Euclidean rank, $rank_{qE}$, one has
\begin{displaymath}
rank_{qE}\Big( X,d_{\Phi}\Big) \; = \; \sum\limits_{i=1}^n \; rank_{qE} \Big( X_i,d_i\Big) .
\end{displaymath} 
\end{corollary}

{\bf Acknowledgement:} It is a pleasure to thank Janko Latschev for useful discussions. 

\section{Non Standard Metric Products}

\label{non-standard}

On $Q^n$ we define a partial ordering $\le$ in the following way: if 
$q^1=(q_1^1,...,q_n^1)$ and $q^2=(q_1^2,...,q_n^2)$ then
\begin{displaymath}
q^1 \; \le \; q^2 \;\;\; :\Longleftrightarrow \;\;\; q_i^1 \; \le q_i^2 \;\; 
\forall i\in \{ 1,2,...,n\} .
\end{displaymath}
Let $\Phi : Q^n \longrightarrow [0,\infty )$ be a function and consider the function
$d_{\Phi}:X\times X \longrightarrow [0,\infty )$,
\begin{displaymath}
d_{\Phi}\Big( (x_1,...,x_n),(y_1,...,y_n)\Big) \; = \; 
\Phi \Big( d_1(x_1,y_1),...,d_n(x_n,y_n)\Big) .
\end{displaymath} 
In order that  $d_{\Phi}$ will be a metric, we clearly have to assume
\begin{description}
\item[(A)] $\Phi (q) \; \ge \; 0 \;\; \forall q\in Q\;\;\;$ and $\;\;\;\Phi (q) \; = \; 0 \;\; 
\Leftrightarrow \;\; q=0$.
\end{description}
The symmetry of $d_{\Phi}$ is obvious. We now translate the triangle inequality for
$d_{\Phi}$ into a condition on $\Phi$. \\
Let $x=(x_1,...,x_n)$, $y=(y_1,...,y_n)$, $z=(z_1,...,z_n)\in X$ and consider the 
``distance vectors''
\begin{eqnarray*}
q^1 & := & \Big( d_1(x_1,z_1),...,d_n(x_n,z_n)\Big), \\
q^2 & := & \Big( d_1(x_1,y_1),...,d_n(x_n,y_n)\Big) \hspace{1cm} \mbox{and} \\
q^3 & := & \Big( d_1(y_1,z_1),...,d_n(y_n,z_n)\Big)
\end{eqnarray*}
in $Q^n$. Since for every $i\in \{ 1,...,n\}$, $x_i,y_i,z_i$ are points in $X_i$
we see that $q^j\le q^k+q^l$ for every permutation $\{ j,k,l \}$ of $\{ 1,2,3 \}$. \\
Now $d_{\Phi}$ satisfies the triangle inequality if $\Phi$ satisfies
\begin{description}
\item[(B)] for all points $q^1,q^2,q^3\in Q^n$ with $q^j\le q^k+q^l$ we have
\begin{displaymath}
\Phi (q^j) \; \le \;  \Phi (q^k) \; + \; \Phi (q^l).
\end{displaymath}
\end{description}

{\bf \underline{Remark}:}
\begin{description}
\item[i)] Note that for $q^1,q^2,q^3$ one can always take a triple of the form 
$p,q,p+q$, hence (B) implies in particular $\Phi (p+q) \le \Phi (p)+ \Phi (q)$, which will be called the sub-additivity of $\Phi$ in the following.
\item[ii)] The condition (B) can be applied for the triple $p,q,q$ in the case
that $p\le 2q$. Then $\Phi (p)\le 2 \Phi (q)$. \\
\end{description}

It is now easy to prove the following result
\begin{lemma} \label{metric}
Let $\Phi : Q^n\longrightarrow [0,\infty )$ be a function. Then $d_{\Phi}$ is a metric
on $X$ for all possible choices of metric spaces $(X_i,d_i)$, $i=1,...,n$, if 
and only if $\Phi$ satisfies (A) and (B).
\end{lemma}
This Lemma still allows strange metrics on a product (even the trivial
product $n=1$). Let for example $\Phi : Q^n\longrightarrow [0,\infty )$ be an 
arbitrary function with $\Phi (0)=0$ and $\Phi (q)\in\{ 1,2\}$,
$\forall q\in Q^n\setminus \{ 0\}$. Then $d_{\Phi}$ is a metric. \\

If we, however, require for example that the product metric space $X$ is always an inner metric space
in the case the $X_i$ are, the conditions on $\Phi$ are very rigid. \\
For the convenience of the reader we recall the notion of an inner metric space.
Let $(X,d)$ be a metric space. For a continuous path $c:[0,1]\longrightarrow X$ 
one defines as usual the length
\begin{equation} \label{arclength}
L(c) \; := \; \sup \Big\{ \sum\limits_{j=1}^k d\Big( c(t_{j-1}),c(t_j) \Big) \Big\},
\end{equation}  
where the $sup$ is taken over all subdivisions
\begin{displaymath}
0 \; = \; t_0 \; \le \; t_1 \; \le ... \le \; t_k \; = \; 1 \hspace{1cm} 
\mbox{of} \;\; [0,1].
\end{displaymath}
$(X,d)$ is called an inner metric space if for all $x,y\in X$, $d(x,y)=\inf L(c)$,
where the $\inf$ is taken over all paths from $x$ to $y$. The curve $c$ is called rectifiable if
L(c) is finite.\\

Spaces of curvature bounded from above (resp. below) are defined via comparison triangles in appropriate space forms of 
constant curvature. We refer the reader to \cite{bh} and \cite{bgp} for details. A space of non-positive Busemann 
curvature is a geodesic space such that the metric is locally convex (\cite{bh}). If the metric is even globally convex, the geodesic space is called convex. For instance each metric space of 
curvature $\leq 0$ has non-positive Busemann curvature and each geodesic CAT($0$)-space is convex (but not vice-verse). 

\vspace{1cm}

We further need the following two Lemmata:
\begin{lemma} \label{psi-norm}
For $\Phi :Q^n\longrightarrow [0,\infty )$ the function
$\Psi : \mathbb{R}^n \longrightarrow [0,\infty )$ defined via
\begin{displaymath}
\Psi \Big( \sum\limits_{i=1}^n \; x_i \, e_i \Big) \; := \; 
\Phi \Big( \sum\limits_{i=1}^n \; |x_i| \, e_i \Big) 
\end{displaymath}
is a norm on $\mathbb{R}^n$ if and only if $\Phi$ satisfies the following conditions:
\begin{description}
\item[(1)] $\Phi (q)\ge 0$ $\forall q\in Q^n$ and $\Phi (q)=0 \; \Leftrightarrow q=0$,
\item[(2)] $\Phi$ is monoton, i.e. $q\le p \; \Longrightarrow \Phi (q) \le \Phi (p)$
$\forall p,q\in Q^n$,
\item[(3)] $\Phi (p+q) \; \le \; \Phi (p) \; + \; \Phi (q)$,
\item[(4)] $\Phi (\lambda q) \; = \; \lambda \Phi (q)$ $\forall p\in Q^n, \lambda \ge 0$.
\end{description}
\end{lemma}
{\bf Proof of Lemma \ref{psi-norm}:} \\
``$\Longrightarrow$'' \\
Let $\Phi$ satisfy $(1)-(4)$. Then
$\Psi \ge 0$, $\Psi (x)=0 \; \Longleftrightarrow \;x=0$ and 
$\Psi (\lambda x) \; = \; |\lambda | \Psi (x)$ directly follow from the 
definition of $\Psi$. In order to verify the subadditivity, note that
for $x=(x_1,...,x_n)$ and $y=(y_1,...,y_n)$
\begin{eqnarray*}
\Psi (x \; + \; y) & = & \Phi \Big( \sum\limits_{i=1}^n \; |x_i+y_i| \, e_i \Big) \\
& \stackrel{(2)}{\le} & \Phi \Big( \sum\limits_{i=1}^n \; (|x_i|+|y_i|) \, e_i \Big) \\
& \stackrel{(3)}{\le} & \Phi \Big( \sum\limits_{i=1}^n \; |x_i| \, e_i \Big) \; + \; 
\Phi \Big( \sum\limits_{i=1}^n \; |y_i| \, e_i \Big) \\
& =  & \Psi (x) \; + \; \Psi (y).
\end{eqnarray*}
``$\Longleftarrow$'' \\
Assume now that $\Psi$ is a norm. Then $\Phi$ clearly satisfies $(1),(3),(4)$. To prove
(2) it is enough to show that $\Phi (p+\lambda e_i)\ge \Phi (p)$ for any unit vector 
$e_i$ and $\lambda \ge 0$. Assume that $\Phi (p+\lambda e_i)< \Phi (p)$. Write
$p=(p_1,...,p_n)\in Q^n$, let 
$q=(p_1,...,p_{i-1},-p_i-\lambda ,p_{i+1},...,p_n)\in \mathbb{R}^n$. Then
$\Psi (q)=\Psi (p+\lambda e_i)< \Psi (p)$ but $p$ is on the segment between $q$ and 
$p+\lambda e_i$. This contradicts the subadditivity of $\Psi$. 
\begin{flushright}
{\bf q.e.d.}
\end{flushright}

\begin{lemma} \label{psi-sp}
Let $\Phi : Q^n \longrightarrow [0,\infty )$ satisfy $(1)-(4)$ as in Lemma \ref{psi-norm}.
Then $\Psi$ as in Lemma \ref{psi-norm} is induced by a scalar product $g_{\Psi}$
on $\mathbb{R}^n$ if and only if $\Phi$ satisfies the property
\begin{displaymath}
(5) \hspace{1cm}
{\Phi}^2 \Big( \sum\limits_{i=1}^n \; {\lambda}_i \, e_i\Big) \; = \; 
\sum\limits_{i=1}^{n} \; {\Phi}^2 ({\lambda}_i \, e_i) \hspace{1cm} \forall {\lambda}_i >0.
\end{displaymath}
In this case the set $\{ e_1,...,e_n \}$ is an orthogonal system of $g_{\Psi}$.
\end{lemma}
{\bf  Proof of Lemma \ref{psi-sp}:} \\
From Lemma \ref{psi-norm} we know that $\Psi$ is a norm if and only if
conditions $(1)-(4)$ hold. Now we show that $\Psi$ satisfies the parallelogram
equation if and only if condition $(5)$ also holds: \\
``$\Longrightarrow$'' \\
Suppose that $\Phi$ satisfies condition $(5)$. Then for $x=(x_1,...,x_n)$ and 
$y=(y_1,...,y_n)$ the parallelogram equation is equivalent to
\begin{displaymath}
\sum\limits_{i=1}^n \; \Big[ |x_i+y_i|^2 \; + \; |x_i-y_i|^2 \; - \; 
2 \, \big[ |x_i|^2 \, + \, |y_i|^2 \big] \Big] {\Phi}^2 (e_i) \; = \; 0 ,
\end{displaymath} 
which holds trivially. \\
``$\Longleftarrow$'' \\
Now suppose that the parallelogram equation holds. For $x=(x_1,...,x_{n-1},0)$ and
$y=(0,...,0,y_n)$ it takes the form 
\begin{displaymath}
{\Phi}^2 \Big( \sum\limits_{i=1}^{n-1} \; |x_i| \, e_i \; + \; |y_n| \, e_n \Big) \; = \;
{\Phi}^2 \Big( \sum\limits_{i=1}^{n-1} \; |x_i| \, e_i \Big) \; + \; 
{\Phi}^2 \Big( |y_n| \, e_n \Big) .
\end{displaymath}
The same computation for $x=(x_1,...,x_{n-2},0,0)$, $y=(0,...,0,y_{n-1},0)$ and so on finally 
yields condition $(5)$.
\begin{flushright}
{\bf q.e.d.}
\end{flushright}

From Lemmata \ref{psi-norm} and \ref{psi-sp} we easily conclude the following propositions:
\begin{proposition} \label{prop-norm}
Let $(V_i,||\cdot ||_i)$, $i=1,...,k$ be normed vector spaces and 
$\Phi : Q^n\longrightarrow [0,\infty )$ be a function. Define the function
$||\cdot ||_{\Phi}:V=V_1\times ... \times V_k \longrightarrow [0,\infty )$ through
\begin{displaymath}
\Big| \Big| (v_1,...,v_k)\Big| {\Big|}_{\Phi} \; := \;  \Phi 
\Big( \sum\limits_{i=1}^k \; ||v_i||_i \, e_i \Big) .
\end{displaymath} 
Then $(V,||\cdot ||_{\Phi})$ is a normed vector space for all possible choices of normed 
vector spaces $(V_i,||\cdot ||_i)$ if and only if $\Psi$ as defined in Lemma 
\ref{psi-norm} is a norm.
\end{proposition}
and
\begin{proposition} \label{prop-sp}
Let $(V_i,||\cdot ||_i)$, $i=1,...,k$, be normed vector spaces the norms of which 
are induced by scalar products $<\cdot ,\cdot >_i$ on $V_i$ and
$\Phi : Q^n\longrightarrow [0,\infty )$ be a function. \\
Then the norm $||\cdot ||_{\Phi}$ on $V=V_1\times ... \times V_k$ as in Proposition \ref{prop-norm} is 
induced by a scalar product $<\cdot ,\cdot >_{\Phi}$ for all choices
of vector spaces $V_i$ with scalar products $<\cdot ,\cdot >_i$, if and only if
the norm $\Psi$ as defined in Lemma \ref{psi-norm} is induced by a scalar product
$g_{\Psi}$ on $\mathbb{R}^n$. \\
Thus for two vectors $v=(v_1,...,v_k),w=(w_1,...w_k)\in v$ one always has
\begin{displaymath}
<v,w>_{\Phi} \; = \; \sum\limits_{i=1}^k \; {\Phi}^2 (e_i) \; <v_i,w_i> \; ,
\end{displaymath}
which is the usual Euclidean product up to a scale of the scalar products on the factors.
\end{proposition}
Note that the degree to that $\{ e_1,...,e_n\}$ fails to be an orthonormal basis
of $g_{\Psi}$ is the degree to that $<\cdot ,\cdot >_{\Phi}$ differs from the
standard scalar product of Euclidean products. \\
{\bf  Outline of the proofs of Propositions \ref{prop-norm} and \ref{prop-sp}:} \\
The ``if parts'' are obvious consequences of the Lemmata \ref{psi-norm} and \ref{psi-sp}. For the only if part one 
considers 
special settings in which all of the factors are the reals with the standard norm (scalar product, respectively). While
conditions $(1)$ and $(4)$ from Lemma \ref{psi-norm} hold trivially, the conditions $(2)$ and $(3)$ now easily follow by
simple constructions in $\mathbb{R}$. Finally, for Theorem \ref{prop-sp}, condition $(5)$ follows just as 
in Lemma \ref{psi-sp}. \\

We are now able to give the 

{\bf Proof of Theorem \ref{theo-character}:} \\
``$\Longrightarrow$'' 
For the ``only if'' part we will restrict our attention to appropriately chosen factors $(X_i,d_i)$.
Of course we do so by assuming that all the factors are the reals with the standard metric: $(X_i,d_i)=(\mathbb{R} 
,d_e)$,
$i=1,...,k$. \\

Similar as in Lemma \ref{metric} we see that $\Phi$ satisfies $(A)=(1)$ and $(B)$ which implies the
subadditivity of $\Phi$.\\

We first show that in all the cases i) - v) $\Phi$ must be induced by a norm. To start with, let us show that $\Phi$ has to be continous. Suppose it is not. By subadditivity, one easily gets that the restriction of $\Phi$ to one of the coordinate lines (say the first one) is discontinous at $0$. Then any two points with different projections on $X_1$ can not be joined by a continous curve in the product space, which therefore can not be a length space.     

That $\Phi$ is induced by a norm now follows immediately from Theorem 6 in \cite{be}, which implies 
that any homogeneous length
metric on $\mathbb{R}^n$ that induces the same topology as the standard metric must come from a norm. \\

$v)$ now follows from the fact that a normed vector space has a lower or upper curvature bound 
if and only if it is Euclidean.
Thus, due to Proposition \ref{prop-sp}, $\Phi$ must be induced by a scalar product. \\

In order to complete the proof in the cases iii) and iv) it suffices to remark that if $\Phi$ does not admit a strictly 
convex unit ball, then the product space is generally not uniquely geodesic. This can be seen for instance in the case  
$(\mathbb{R}^n,\Psi )$.  \\

``$\Longleftarrow$ '' 
Let now $\Phi$ be induced by a norm. Note that from Lemma \ref{psi-norm} it 
follows immediately that $\Phi$ satisfies the conditions $(1)-(4)$ of Lemma \ref{psi-norm}. 
In order to show that for any choices of 
inner metric (geodesic) spaces $X_1,...,X_k$ the product $(X,d_{\Phi})$ is an inner metric (geodesic) space we 
prove the following
\begin{lemma} \label{productlength}
Let $(X_i,d_i)$ be metric spaces and $c_i: [0,1]\longrightarrow X_i$ be continuous curves
parameterized by arclength connecting $p_i\in X_i$ with 
$q_i\in X_i$, $i=1,...,k$. Denote by $l_i$ the $(X_i,d_i)$-length of $c_i$ and suppose that
$\Phi$ satisfies conditions $(1)-(4)$. Then the $(X,d_{\Phi})$-length of the product 
curve $c=(c_1,...,c_k):[0,1]\longrightarrow X$ is $L(c)=\Phi(l_1,...,l_k)$. \\
Furthermore $c$ is also parameterized by arclength.
\end{lemma}
{\bf Proof of Lemma \ref{productlength}:} \\
Note that the $(X_i,d_i)$-length $L_{c_i}$ of $c_i$ is given through
\begin{displaymath}
l_i \; = \; L(c_i) \; = \; \lim\limits_{N\longrightarrow \infty}\sum\limits_{j=1}^N \;
d_i \Big( c_i(\frac{j-1}{N}), c_i(\frac{j}{N})\Big) ,
\end{displaymath}
where $d_i\big( c_i(\frac{j-1}{N}), c_i(\frac{j}{N})\big) \leq \frac{l_i}{N}$. \\
For the $(X,d_{\Phi})$-length $L(c)$ of $c$ one has
\begin{eqnarray*}
L(c) & = & \lim\limits_{N\longrightarrow \infty}\sum\limits_{j=1}^N \;
d_{\Phi} \Big( c(\frac{j-1}{N}), c(\frac{j}{N})\Big) \\
& = &\lim\limits_{N\longrightarrow \infty}\sum\limits_{j=1}^N \; 
\Phi \Big( d_1\Big( c_1(\frac{j-1}{N}),c_1(\frac{j}{N})\Big) ,...,
d_k\Big( c_k(\frac{j-1}{N}),c_k(\frac{j}{N})\Big) \Big) \\
& \stackrel{(2)}{\le} & \lim\limits_{N\longrightarrow \infty}\sum\limits_{j=1}^N \; 
\Phi \Big( \frac{l_1}{N},...,\frac{l_k}{N}\Big) \\
& \stackrel{(4)}{=} & \lim\limits_{N\longrightarrow \infty}\sum\limits_{j=1}^N \; 
\frac{1}{N} \Phi (l_1,...,l_k) \\
& = & \Phi (l_1,...,l_k).
\end{eqnarray*}
On the other hand the continuity and subadditivity of $\Phi$ yield:
\begin{eqnarray*}
L(c) & = & \lim\limits_{N\longrightarrow \infty}\sum\limits_{j=1}^N \;
d_{\Phi} \Big( c(\frac{j-1}{N}), c(\frac{j}{N})\Big) \\
& = &\lim\limits_{N\longrightarrow \infty}\sum\limits_{j=1}^N \; 
\Phi \Big( d_1\Big( c_1(\frac{j-1}{N}),c_1(\frac{j}{N})\Big) ,...,
d_k\Big( c_k(\frac{j-1}{N}),c_k(\frac{j}{N})\Big) \Big) \\
& \stackrel{(3)}{\ge} & \lim\limits_{N\longrightarrow \infty} 
\Phi \Big( \sum\limits_{j=1}^N \; d_1\Big( c_1(\frac{j-1}{N}),c_1(\frac{j}{N})\Big) ,...,
\sum\limits_{j=1}^N \; d_k\Big( c_k(\frac{j-1}{N}),c_k(\frac{j}{N})\Big) \Big) \\
& = & \Phi (l_1,...,l_k),
\end{eqnarray*}
where the last equality is due to the continuity of $\Phi$.
\begin{flushright}
$\Box$
\end{flushright}
Let now $(X_i,d_i)$, $i=1,...,k$, be inner metric spaces. Then the distance of any two points
$p_i,q_i\in (X_i,d_i)$ may be approximated arbitrarily good by the lengths of continuous 
curves in $(X_i,d_i)$ joining $p_i$ and $q_i$. Thus $(X,d_{\Phi})$ turns out to be
an inner metric space itself, due to the definition of $d_{\Phi}$, the validity of Lemma
\ref{productlength} and the continuity of $\Phi$. \\
This completes the proof in the cases i) and ii). \\

v) now just follows from Proposition \ref{prop-sp} and the fact that stretching the factor metrics 
$d_i$ by ${\Phi}(e_i)$ and then taking the standard product yields as product a space with an upper (resp. lower) curvature bound. 
This last fact can be seen by an easy comparison argument using the fact that the Euclidean product of two constant 
curvature space forms has upper and lower curvature bounds. See also \cite{bh}, II 1.16. and \cite{bgp}. \\

The local convexity of the distance function in iv) is obvious. 
In order to finish the proof in the cases iii) and iv) we have to show that the product is uniquely geodesic again and 
we might as well only consider the product of two factors in order to avoid unneccesary index complications: \\
Let therefore $(x_1,x_2), (y_1,y_2) \in X$, set $D:=d_{\Phi}((x_1,x_2),(y_1,y_2))$ and let 
$c=(c_1,c_2):[0,D] \longrightarrow X$ be a geodesic in $X$ joining $(x_1,x_2)$ to
$(y_1,y_2)$. 
Set $v_1:=(d_1(x_1,c_1(t)),d_2(x_2,c_2(t))) \in Q^2$, $v_2:=(d_1(y_1,c_1(t)),d_2(y_2,c_2(t))) \in Q^2$ and 
$v:=(d_1(x_1,y_1),d_2(x_2,y_2)) \in Q^2$. By the triangle inequality in each component, we have $v_1+v_2 \geq v$. With 
$(c_1,c_2)$ being a geodesic and by subadditivity, we must have 
\begin{displaymath}
\Phi(v_1)+\Phi(v_2)=\Phi(v)\leq \Phi(v_1+v_2) \leq \Phi(v_1)+\Phi(v_2).
\end{displaymath}
Since the norm ball associated to $\Phi$ is strictly convex, the conditions $\Phi(v_1)+\Phi(v_2)=\Phi(v)$ and $v\leq 
v_1+v_2$ can only be satisfied if $v_1=\lambda_1 v, v_2=\lambda_2 v$ with $\lambda_1+\lambda_2=1$. Since $\Phi(v_1)=t, 
\Phi(v)=D$, we get $\lambda_1=\frac{t}{D},\lambda_2=1-\frac{t}{D}$. The spaces $X_1,X_2$ are uniquely geodesic and 
therefore $c_1(t),c_2(t)$ are fixed by these equations. Hence there is a unique geodesic $(c_1,c_2)$ joining $(x_1,x_2)$ 
and $(y_1,y_2)$. 
\begin{flushright}
{\bf q.e.d.}
\end{flushright}


\section{Minkowski Rank of Products}

\label{sec-minkowski}

In this section we prove Theorem \ref{generalized-minkadd} and Corollary \ref{quasiadd}. As the proof of Theorem \ref{generalized-minkadd}
is almost the same as the one of Theorem 2 of \cite{foesch} we keep it fairly short. \\

In \cite{foesch} we introduced the notions of the Euclidean and the Minkowski rank for arbitrary metric spaces as follows.

\begin{definition} Minkowski- and Euclidean rank for metric spaces
\begin{description}
\item[a)] For an arbitrary metric space $(X,d)$ the  {\bf Minkowski rank} is
\begin{displaymath}
rank_{M}(X,d) \; := \; 
\sup_{(V,||\cdot ||)} \Big\{ dim V \; \Big| \; \exists \; \mbox{isometric map} \;\; 
i_V:(V,||\cdot ||) \longrightarrow (X,d) \Big\} .
\end{displaymath}
\item[b)] The {\bf Euclidean rank} is defined as
\begin{displaymath}
rank_{E}(X,d) \; := \; 
\sup \Big\{ n \in \mathbb{N} \; \Big| \; \exists \; \mbox{isometric map} \;\; 
i_{\mathbb{E}^n}:\mathbb{E}^n\longrightarrow (X,d) \Big\} .
\end{displaymath}
\end{description}
\end{definition} 

The Minkowski rank was shown to be additive with respect to the standard product of arbitrary metric spaces, whereas a counterexample to the corresponding additivity of the Euclidean rank was provided. Since for metric spaces of locally one-side bounded Alexandrov curvature these two ranks coincide (\cite{ri},\cite{foesch}), the additivity of the Euclidean rank with respect to the standard product follows for instance for metric spaces of non-positive or non-negative Alexandrov curvature. \\

In \cite{kl} Kleiner considered another notion of rank which we will refer to as quasi-Euclidean rank in the following. Recall that given two metric spaces $(X,d_X),(Y,d_Y)$, a map $f:X \longrightarrow Y$ is called quasi-isometric embedding of $X$ in $Y$, if there exist $\lambda\geq 1$ and $\epsilon>0$ such that for all $x_1,x_2 \in X$:
\begin{displaymath}
\lambda^{-1} d_X(x_1,x_2) -\epsilon \leq d_Y(f(x_1),f(x_2)) \leq \lambda d_X(x_1,x_2) +\epsilon
\end{displaymath}

\begin{definition} 
 The {\bf quasi-Euclidean rank} is defined as
\begin{displaymath}
rank_{qE}(X,d) \; := \; 
\sup \Big\{ n \in \mathbb{N} \; \Big| \; \exists \; \mbox{quasi-isometric embedding} \;\; 
f_{\mathbb{E}^n}:\mathbb{E}^n\longrightarrow (X,d) \Big\} .
\end{displaymath}
\end{definition} 
 
{\bf Scetch of proof of Theorem \ref{generalized-minkadd}:}\\
The analogue of Proposition 2 in \cite{foesch} is
\begin{proposition} \label{generalized-prop}
Let $A$ denote an affine space on which the normed vector space $(V,|\cdot |)$ acts
simply transitively. Let further $(X_i,d_i)$, $i=1,...,n$, be metric spaces, 
$\Phi :Q^n \longrightarrow [0,\infty )$ be a function satisfying
conditions $(1)-(4)$ such that the norm ball of $\Psi$ is strictly convex and 
let $\varphi : (A,|\cdot |) \longrightarrow ({\Pi}_{i=1}^nX_i,d_{\Phi})$ be an 
isometric map. Then there exist pseudonorms $||\cdot ||_i$, $i=1,...,n$, on V such
that
\begin{description}
\item[i)] $|v| \; = \; \Phi \Big( \sum\limits_{i=1}^n \; ||v||_i \, e_i \Big)
\hspace{1cm} \forall v\in V \hspace{1cm} \mbox{and}$
\item[ii)] ${\varphi}_i:(A,||\cdot ||_i)\longrightarrow (X_i,d_i)$, $i=1,...,n$ are
isometric.
\end{description}
\end{proposition} 

Note that from Proposition \ref{generalized-prop} one achieves the subadditivity of the Minkowski rank, while
the superadditivity just follows by considering the product of two isometric embeddings of normed 
vector spaces into the factors. \\
We define ${\alpha}_i:A\times V \longrightarrow [0,\infty)$, $i=1,...,n$, via
\begin{displaymath}
{\alpha}_i(a,v) \; := \; d_i\Big( {\varphi}_i(a),{\varphi}_i(a+v) \Big).
\end{displaymath}
Since $\varphi$ is isometric we have
\begin{equation} \label{lemma1'1}
\Phi \Big( \sum\limits_{i=1}^n \; {\alpha}_i(a,v) \, e_i \Big) \; := \; 
d_{\Phi} \Big( \varphi (a),\varphi (a+v)\Big) \; = \; |v|.
\end{equation}
In order to prove Proposition \ref{generalized-prop}, we might as well restrict to the case $n=2$ and show that 
$||\cdot ||_i:V\longrightarrow \mathbb{R}^+$,
$||v||_i:={\alpha}_i(v)$ $\forall v\in V$, $i=1,2$, are pseudonorms on $V$. Therefore we establish the following 
properties of $\alpha$:
\begin{description}
\item[1)]
\begin{displaymath}
{\alpha}_i(a,v) \; = \; {\alpha}_i(a+v,v), \;\; i=1,2, \hspace{2cm} \forall a\in A, v\in V,
\end{displaymath}
\item[2)]
\begin{displaymath}
{\alpha}_i(a,tv) \; = \; |t| {\alpha}_i(a,v), \;\; i=1,2, \hspace{2cm} \forall a\in A, v\in V, t\in \mathbb{R},
\end{displaymath} 
\item[3)]
\begin{displaymath}
{\alpha}_i(a,v) \; = \; {\alpha}_i(b,v), \;\; i=1,2, \hspace{2cm} \forall a,b\in A, v\in V \;\; \mbox{and}
\end{displaymath}
\item[4)]
\begin{displaymath}
{\alpha}_i(v+w) \; \le \; {\alpha}_i(v) \; + {\alpha}_i(w), \;\; i=1,2, \hspace{2cm} \forall v,w\in V,
\end{displaymath}
where ${\alpha}_i(v):={\alpha}_i(a,v)$ with $a\in A$ arbitrary (compare with $3)$). 
\end{description}
In order to prove $1)$ we note that the $d_i$'s triangle inequality yields
\begin{displaymath}
{\alpha}_i(a,v) \; + \; {\alpha}_i(a+v,v) \; \ge \; {\alpha}_i(a,2v). 
\end{displaymath}
Therefore the monotonicity of $\Phi$ gives
\begin{equation} \label{lemma1'2}
\Phi \Big( \sum\limits_{i=1}^n \; [{\alpha}_i(a,v) \; + \; {\alpha}_i(a+v,v) ] \, e_i \Big)
\; \ge \; 
\Phi \Big( \sum\limits_{i=1}^n \; {\alpha}_i(a,2v) \, e_i \Big)
\; = \; 
|2v| \; = \; 2 \, |v|.
\end{equation}
Set $x:=\sum_{i=1}^n{\alpha}_i(a,v)e_i$ and $y:=\sum_{i=1}^n{\alpha}_i(a+v,v)e_i$ and
note that with equations (\ref{lemma1'1}) and (\ref{lemma1'2}) one has
\begin{displaymath}
\Phi (x) \; = \; |v| \; = \; \Phi (y) \hspace{1cm} \mbox{and} \hspace{1cm}
\Phi (x+y) \; = \; 2 \; |v| 
\end{displaymath}
and hence
\begin{displaymath}
\Phi (x+y) \; = \; \Phi (x) \; + \; \Phi (y) \; = \; 2\Phi (x) \; = \; 2\Phi (y).
\end{displaymath}
From this it follows with the strict convexity of $\Phi$ that
\begin{displaymath}
x \; = \; y \hspace{0.5cm} \Longleftrightarrow \hspace{0.5cm} 
{\alpha}_i(a,v) \; = \; {\alpha}_i(a+v,v) \;\;\; \forall i=1,...,n,
\end{displaymath}
which proves $1)$. \\
In order to prove $2)$ we note that the $d_i$'s triangle inequality yields for all $n\in \mathbb{N}$
\begin{displaymath}
{\alpha}_i(a,nv) \; \le \; \sum\limits_{k=0}^{n-1}{\alpha}_i(a+kv,v) \; = \; n{\alpha}_i (a,v),
\end{displaymath}
where the last equation follows from $1)$ by induction. Thus we find $\forall n\in \mathbb{N}, v\in V, a\in A$:
\begin{displaymath}
n^2||v||^2 \; = \;
\Phi \Big( {\alpha}_1(a,nv), {\alpha}_2(a,nv) \Big) \; \le \;
n^2 \Phi \Big( {\alpha}_1(a,v), {\alpha}_2(a,v) \Big) 
\; = \; n^2||v||^2 
\end{displaymath}
and therefore
\begin{displaymath}
{\alpha}_i(a,nv) \; = \; n{\alpha}_i (a,v), \;\; i=1,2, \hspace{2cm} \forall n\in \mathbb{N}, v\in V,
a\in A. 
\end{displaymath}
The claim now follows by the usual extension to $n\in \mathbb{Q}$ and finally to $n\in \mathbb{R}$. \\
In order to prove $3)$ we observe that for $n\in \mathbb{N}$ we have
\begin{eqnarray*}
\Big|{\alpha}_i(a,nv) \; - \; {\alpha}_i(b,nv)\Big| & = & 
\Big| d_i\Big( {\varphi}_i(a),{\varphi}(a+nv)\Big) - d_i\Big( {\varphi}_i(b),{\varphi}(b+nv)\Big) \Big| \\
& \le & d_i \Big( {\varphi}_i(a),{\varphi}_i(b)\Big) \; + \;  
d_i \Big( {\varphi}_i(a+nv),{\varphi}_i(b+nv)\Big) \\
& \le &  d \Big( {\varphi}(a),{\varphi}(b)\Big) \; + \; d \Big( {\varphi}(a+nv),{\varphi}(b+nv)\Big) \\
& = & 2||b-a||, \;\;\;\;\;\;\; i=1,2,
\end{eqnarray*}
and therefore
\begin{displaymath}
{\alpha}_i(a,v) \; = \; \lim\limits_{n\longrightarrow \infty} \frac{1}{n}{\alpha}_i(a,nv)
\; = \;  \lim\limits_{n\longrightarrow \infty} \frac{1}{n}{\alpha}_i(b,nv) \; = \; 
{\alpha}_i(b,v), \;\; i=1,2.
\end{displaymath}
Finally $4)$ follows from
\begin{displaymath}
{\alpha}_i(v+w) \; = \; {\alpha}_i(a,v+w) \; \le \; {\alpha}_i(a,v) \; + \; {\alpha}_i(a+v,w) \; = \; 
{\alpha}_i(v) \; + \; {\alpha}_i(w),
\end{displaymath}
where the inequality follows by the $d_i$'s triangle inequality and the last equation is due to
$3)$. \\
 
The following example shows that the strict convexity of $\Phi$, assumed in Theorem
\ref{generalized-minkadd}, is a 
necessary condition. Take the interval $\mathbb{R}_+:=[0,\infty)$ with the metric induced from $\mathbb{R}$. By simple 
geometric arguments, $rank_{M}(\mathbb{R}_+)=0$. However, with $\Phi(x_1,x_2):=x_1+x_2$, the $\Phi$-product of 
two copies of $\mathbb{R}_+$ admits the geodesic line $c(t)=(-t,0)$ if $t\leq 0$, $c(t)=(0,t)$ if $t \geq 0$.\\ 

{\bf Proof of Corollary \ref{quasiadd}:}\\ 
It is evident that the product of locally compact, convex metric spaces with cocompactly acting isometry group also
satisfies these conditions.
From Theorem D of \cite{kl} it follows that under these conditions the quasi-Euclidean and the Minkowski rank coincide. Thus the validity of Corollary \ref{quasiadd} is a consequence of Theorem \ref{generalized-minkadd}. 
\begin{flushright}
{\bf q.e.d.} 
\end{flushright}

Remark: In general, the quasi-Euclidean rank is not additive with respect to the standard product of metric spaces. Just consider the above example of two copies of $\mathbb{R}_+$ with the standard product. The geodesic defined above is then a quasi-geodesic.


{\footnotesize UNIVERSIT\"AT Z\"URICH, MATHEMATISCHES INSTITUT, WINTERTHURERSTRASSE 190, 
CH-8057 Z\"URICH, SWITZERLAND \\
E-mail addresses: bernig,foertsch,vschroed@math.unizh.ch

\end{document}